\documentclass[amssymb,amstex,amsmath,10pt,a4paper]{amsart}
\usepackage{amsmath,amssymb,latexsym,amsfonts,url}

\newtheorem{thm}{Theorem}[section]
\newtheorem*{theorem}{Theorem}

\newtheorem{lem}[thm]{Lemma}

\newcommand{\Ric}{{\rm Ric}}

\numberwithin{equation}{section}

\begin{document}
\title[Rigidity of minimal submanifolds in hyperbolic space]{Rigidity of minimal submanifolds \\in hyperbolic space}
\author{Keomkyo Seo}
\address{Department of Mathematics\\
Sookmyung Women's University\\
Hyochangwongil 52, Yongsan-ku \\
Seoul, 140-742, Korea}
\email{kseo@sookmyung.ac.kr}


\begin{abstract}
We prove that if an $n$-dimensional complete minimal submanifold $M$ in hyperbolic space has sufficiently small total scalar curvature then $M$ has only one end. We also prove
that for such $M$ there exist no nontrivial $L^2$ harmonic $1$-forms on $M$.\\

\end{abstract}

\maketitle

\section{Introduction}
A complete minimal hypersurface $M \subset \overline{M}^{n+1}$ is said to be $stable$ if the second variation of its volume is always nonnegative for any normal variation with compact support. In \cite{CSZ}, Cao, Shen, and Zhu showed that a complete connected stable minimal hypersurface in Euclidean space must have exactly one end. Later Ni\cite{Ni} proved that if an $n$-dimensional complete minimal submanifold $M$ in Euclidean space has sufficiently small total scalar curvature (i.e., $\int_M |A|^n dv < \infty$) then $M$ has only one end. More precisely, he proved
\begin{theorem} {\rm(\cite{Ni})}
Let $M$ be an $n$-dimensional complete immersed minimal
submanifold in $\mathbb{R}^{n+p}$, $n \geq 3$. If
\begin{eqnarray*}
\Big(\int_M |A|^n dv \Big)^{1 \over n} < C_1 = \sqrt{{n \over
{n-1}}C_s^{-1}},
\end{eqnarray*}
then $M$ has only one end. (Here $C_s$ is a Sobolev constant in
\cite{HS}.)
\end{theorem}
Recently the author\cite{{Seo}} improved the upper bound $C_1$ of total scalar curvature in the above theorem.

In this paper, we shall prove that the analogue of the above theorem still holds in hyperbolic space. Throughout this paper, we shall denote by
$\mathbb{H}^n$ the $n$-dimensional hyperbolic space of constant sectional curvature $-1$. Our main result is the following.
\begin{thm} \label{main1}
Let $M$ be an $n$-dimensional complete immersed minimal submanifold in $\mathbb{H}^{n+p}$, $n \geq 5$. If the total scalar curvature satisfies
\begin{eqnarray*}
\Big(\int_M |A|^n dv \Big)^{1 \over n} < {1 \over
{n-1}}\sqrt{n(n-4)C_s^{-1}},
\end{eqnarray*}
then $M$ has only one end.
\end{thm}
Miyaoka \cite{Miyaoka} showed that if $M$ is a complete stable minimal hypersurface in $\mathbb{R}^{n+1}$, then there are no nontrivial $L^2$ harmonic $1$-forms on $M$. In \cite{Yun}, Yun proved that if $M \subset \mathbb{R}^{n+1}$ is a complete minimal hypersurface with $\displaystyle{\Big(\int_M |A|^n dv \Big)^{1 \over n} < C_2 = \sqrt{C_s^{-1}}}$, then there are no nontrivial $L^2$ harmonic $1$-forms on $M$. Recently the author\cite{Seo} showed that this result is still true for any complete minimal submanifold $M^n$ in $\mathbb{R}^{n+p}$, $p \geq 1$. We shall prove that if $M$ is an $n$-dimensional complete minimal submanifold with sufficiently small total scalar curvature in hyperbolic space, then there exist no nontrivial $L^2$ harmonic $1$-forms on $M$. More precisely, we prove
\begin{thm}  \label{main2}
Let $M$ be an $n$-dimensional complete immersed minimal submanifold in $\mathbb{H}^{n+p}$, $n \geq 5$. If
\begin{eqnarray*}
\Big(\int_M |A|^n dv \Big)^{1 \over n} < \sqrt{{n(n-4) \over
(n-1)^2} \frac{1}{C_s}} ,
\end{eqnarray*}
then there are no nontrivial $L^2$ harmonic 1-forms on $M$.
\end{thm}

\section{Proof of the theorems}
To prove Theorem \ref{main1}, we begin with the following useful facts.
\begin{lem}[Sobolev inequality {\rm \cite{HS}}]   \label{MS}
Let $M$ be an $n$-dimensional complete immersed minimal submanifold in
$\mathbb{H}^{n+p}$, $n \geq 3$. Then for any $\phi \in W_0
^{1,2}(M)$ we have
\begin{eqnarray*}
\Big(\int_M |\phi|^{\frac{2n}{n-2}} dv\Big)^{\frac{n-2}{n}} \leq C_s
\int_M |\nabla \phi|^2 dv,
\end{eqnarray*}
where $C_s$ depends only on $n$.
\end{lem}

\begin{lem} {\rm (\cite{Leung})}  \label{Leung}
Let $M$ be an $n$-dimensional complete immersed minimal submanifold in
$\mathbb{H}^{n+p}$. Then the Ricci curvature of $M$ satisfies
\begin{eqnarray*}
\Ric(M) \geq -(n-1) -{n-1 \over n}|A|^2.
\end{eqnarray*}
\end{lem}
\noindent Recall that the first eigenvalue of a Riemannian manifold $M$ is defined as
\begin{align*}
\lambda_1 (M) = \inf_{f} \frac{\int_M |\nabla f|^2}{\int_M f^2} ,
\end{align*}
where the infimum is taken over all compactly supported smooth functions on $M$. For a complete stable minimal
hypersurface $M$ in $\mathbb{H}^{n+1}$, Cheung and Leung
\cite{Cheung and Leung} proved that
\begin{eqnarray} \label{thm : CL}
\frac{1}{4}(n-1)^2 \leq \lambda_1 (M) .
\end{eqnarray}
Here this inequality is sharp because equality holds when $M$ is
totally geodesic (\cite{McKean}).

Let $u$ be a harmonic function on $M$. From Bochner formula, we have
\begin{eqnarray*}
{1 \over 2}\Delta(|\nabla u|^2) = \sum {u_{ij}}^2 + \Ric(\nabla u,
\nabla u) .
\end{eqnarray*}
Then Lemma \ref{Leung} gives
\begin{eqnarray*}
{1 \over 2}\Delta(|\nabla u|^2) \geq \sum {u_{ij}}^2  -(n-1)|\nabla u|^2 -{n-1 \over
n}|A|^2 |\nabla u|^2.
\end{eqnarray*}
Choose the normal coordinates at $p$ such that $u_1(p) =
|\nabla u|(p)$, $u_i(p) = 0$ for $i \geq 2$. Then we have
\begin{eqnarray*}
\nabla_j |\nabla u| = \nabla_j ( \sqrt{\sum u_i^2}) = \frac{\sum
u_i u_{ij} }{|\nabla u|} = u_{1j}  \ .
\end{eqnarray*}
Therefore we obtain $|\nabla | \nabla u||^2 = \sum {u_{1j}}^2$.\\

On the other hand, we have
\begin{eqnarray*}
{1 \over 2}\Delta(|\nabla u|^2) = |\nabla u|\Delta |\nabla u| +
|\nabla|\nabla u||^2 .
\end{eqnarray*}
Therefore it follows
\begin{eqnarray*}
\sum {u_{ij}}^2 -(n-1)|\nabla u|^2 -{n-1 \over n}|A|^2 |\nabla u|^2 \leq |\nabla
u|\Delta |\nabla u| + \sum {u_{1j}}^2.
\end{eqnarray*}
Hence we get
\begin{eqnarray*}
|\nabla u|\Delta |\nabla u| + (n-1)|\nabla u|^2+ {n-1 \over n}|A|^2 |\nabla u|^2
&\geq& \sum {u_{ij}}^2 - \sum {u_{1j}}^2 \\
&\geq& \sum_{i \neq 1}{u_{i1}}^2 + \sum_{i \neq 1} {u_{ii}}^2 \\
&\geq&  \sum_{i \neq 1}{u_{i1}}^2 + {1 \over {n-1}}(\sum_{i \neq
1}u_{ii})^2 \\
&\geq& {1 \over {n-1}}{\sum_{i \neq 1} u_{i1}}^2 =  {1 \over
{n-1}}|\nabla|\nabla u||^2 ,
\end{eqnarray*}
where we used $\Delta u = \sum u_{ii} = 0$ in the last inequality.
Therefore we get
\begin{eqnarray} \label{harmonic}
|\nabla u|\Delta |\nabla u| + {n-1 \over n}|A|^2 |\nabla u|^2 + (n-1)|\nabla u|^2\geq {1
\over {n-1}}|\nabla|\nabla u||^2 .
\end{eqnarray}
Now we are ready to prove Theorem \ref{main1}.\\

{\bf Proof of Theorem \ref{main1}.} Suppose that $M$ has at least
two ends. First we note that if $M$ has more than one end then
there exists a nontrivial bounded harmonic function $u(x)$ on $M$
which has finite total energy(\cite{Wei}). Let
$f=|\nabla u|$. From (\ref{harmonic}) we have
\begin{eqnarray*}
f\Delta f + {n-1 \over n } |A|^2 f^2 + (n-1)f^2 \geq {1 \over n-1}|\nabla
f|^2.
\end{eqnarray*}
Fix a point $p \in M$ and for $R>0$ choose a cut-off function
satisfying $0 \leq \varphi \leq 1$, $\varphi \equiv 1$ on
$B_p(R)$, $\varphi = 0 $ on $M \setminus B_p(2R)$, and
$\displaystyle{|\nabla \varphi| \leq {1 \over R}}$, where $B_p(R)$ denotes the ball of radius $R$ centered at $p\in M$. Multiplying
both sides by $\varphi^2$ and integrating over $M$, we have
\begin{eqnarray*}
\int_M \varphi^2 f\Delta f dv + {n-1 \over n } \int_M \varphi^2
|A|^2 f^2 dv + (n-1)  \int_M \varphi^2 f^2 dv \geq {1 \over n-1}\int_M \varphi^2|\nabla f|^2 dv.
\end{eqnarray*}
From the inequality (\ref{thm : CL}), we see
\begin{eqnarray*}
\frac{(n-1)^2}{4} \leq \lambda_1 (M) \leq \frac{\int_M |\nabla (\varphi f)|^2 dv}{\int_M \varphi^2 f^2 dv} .
\end{eqnarray*}
Therefore we get
\begin{eqnarray*}
\int_M \varphi^2 f\Delta f dv + {n-1 \over n } \int_M \varphi^2
|A|^2 f^2 dv + \frac{4}{n-1} \int_M |\nabla (\varphi f)|^2 dv \geq {1 \over n-1}\int_M \varphi^2|\nabla f|^2 dv.
\end{eqnarray*}
Using integration by parts, we get
\begin{align*}
-\int_M |\nabla f|^2 \varphi^2 dv - 2 \int_M f \varphi \langle
\nabla f,\nabla \varphi \rangle dv + {n-1 \over n }\int_M
\varphi^2 |A|^2 f^2 dv + \frac{4}{n-1} \int_M |\nabla (\varphi f)|^2 dv \\
\geq {1 \over n-1}\int_M \varphi^2|\nabla
f|^2 dv.
\end{align*}
Applying Schwarz inequality, for any positive number $a>0$, we
obtain
\begin{align} \label{Sch1}
{n-1 \over n }\int_M \varphi^2 |A|^2 f^2 dv &+ \Big(\frac{4}{n-1}+{1
\over a} + \frac{4}{a(n-1)}\Big) \int_M f^2 |\nabla \varphi|^2 dv  \nonumber \\ &\geq \Big({n \over
{n-1}} - a -\frac{4}{n-1} - \frac{4a}{n-1} \Big)\int_M \varphi^2|\nabla f|^2 dv.
\end{align}
On the other hand, applying Sobolev inequality(Lemma \ref{MS}), we
have
\begin{eqnarray*}
\int_M |\nabla (f\varphi)|^2 dv \geq {C_s}^{-1} \Big( \int_M
(f\varphi)^{\frac{2n}{n-2}} dv \Big)^{n-2 \over n}.
\end{eqnarray*}
Thus applying Schwarz inequality again, we have for any positive
number $b>0,$
\begin{eqnarray} \label{Sch2}
(1+b)\int_M\varphi^2|\nabla f|^2 dv \geq {C_s}^{-1} \Big( \int_M
(f\varphi)^{\frac{2n}{n-2}} dv \Big)^{n-2 \over n} \\ - (1 + {1
\over b})\int_M f^2 |\nabla \varphi|^2 dv.  \nonumber
\end{eqnarray}
Combining (\ref{Sch1}) and (\ref{Sch2}), we have
\begin{align*}
{n-1 \over n }\int_M \varphi^2 |A|^2 f^2 dv  &\geq \Big(\frac{n-4-4a}{n-1} -a \Big)\int_M \varphi^2|\nabla f|^2 dv \\
&\geq  \frac{{C_s}^{-1}}{b+1}\Big(\frac{n-4-4a}{n-1} -a \Big)\Big( \int_M
(f\varphi)^{\frac{2n}{n-2}} dv \Big)^{n-2 \over n} \\
& \quad - \Big\{\frac{1}{b}\Big(\frac{n-4-4a}{n-1} -a \Big)  + \frac{n+3+4a}{a(n-1)}\Big\} \int_M f^2 |\nabla \varphi|^2 dv .
\end{align*}
Applying H\"{o}lder inequality, we get
\begin{eqnarray*}
\int_M \varphi^2 |A|^2 f^2 dv \leq \Big(\int_M |A|^n \Big)^{2
\over n} \Big(\int_M (f\varphi)^{\frac{2n}{n-2}} dv \Big)^{n-2
\over n}.
\end{eqnarray*}
Finally we obtain
\begin{align*}
&\Big\{\frac{1}{b}\Big(\frac{n-4-4a}{n-1} -a \Big) + \frac{n+3+4a}{a(n-1)}\Big\} \int_M f^2 |\nabla \varphi|^2 dv \\
&\geq \Big\{\frac{{C_s}^{-1}}{b+1}\Big(\frac{n-4-4a}{n-1} -a \Big)  - {n-1 \over n }\Big(\int_M |A|^n dv \Big)^{2
\over n}\Big\} \Big(\int_M (f\varphi)^{\frac{2n}{n-2}} dv 
\Big)^{n-2 \over n} .
\end{align*}
By the assumption on the total scalar curvature, we choose $a$ and $b$ small enough such that
\begin{eqnarray*}
\Big\{\frac{{C_s}^{-1}}{b+1}\Big(\frac{n-4-4a}{n-1} -a \Big)  - {n-1 \over n }\Big(\int_M |A|^n dv \Big)^{2
\over n}\Big\} \geq \varepsilon > 0
\end{eqnarray*}
for sufficiently small $\varepsilon>0$.
Then letting $R \rightarrow \infty $, we have $f \equiv 0$, i.e.,
$|\nabla u| \equiv 0$. Therefore $u$ is constant. This contradicts
the assumption that $u$ is a nontrivial harmonic function.   \qed

\medskip

In the proof of Theorem \ref{main1}, if we do not use the fact that $\lambda_1 (M) \geq \frac{(n-1)^2}{4}$ and assume that
\begin{eqnarray*}
\Big(\int_M |A|^n dv \Big)^{2 \over n} < \frac{n}{n-1}{C_s}^{-1}\Big({n \over {n-1}}- \frac{n-1}{\lambda_1 (M)} \Big)
\end{eqnarray*}
for  $\lambda_1 (M) > \frac{(n-1)^2}{n}$, one can see that $M^n$ must have exactly one end by using the same argument as in the above proof, when $n \geq 3$. In other words, it follows

\begin{thm}
Let $M$ be an $n$-dimensional complete immersed minimal submanifold in $\mathbb{H}^{n+p}$, $n \geq 3$. Assume that $\lambda_1 (M) > \frac{(n-1)^2}{n}$ and the total scalar curvature satisfies
\begin{eqnarray*}
\Big(\int_M |A|^n dv \Big)^{2 \over n} < \frac{n}{n-1}{C_s}^{-1}\Big({n \over {n-1}}- \frac{n-1}{\lambda_1 (M)} \Big) .
\end{eqnarray*}
Then $M$ must have only one end.
\end{thm}

\medskip

{\bf Proof of Theorem \ref{main2}.} Let $\omega$ be an $L^2$
harmonic 1-form on minimal submanifold $M$ in $\mathbb{H}^{n+p}$.
Then $\omega$ satisfies
\begin{eqnarray*}
\Delta \omega = 0 \ \ {\rm and}\ \  \int_M |\omega|^2 dv < \infty.
\end{eqnarray*}
It follows from Bochner formula
\begin{eqnarray*}
\Delta|\omega|^2 = 2(|\nabla \omega|^2 + \Ric(\omega,\omega)) .
\end{eqnarray*}
We also have
\begin{eqnarray*}
\Delta|\omega|^2 = 2(|\omega|\Delta |\omega| + |\nabla|\omega||^2)
.
\end{eqnarray*}
Since $\displaystyle{|\nabla \omega|^2 \geq {n \over {n-1}}
|\nabla|\omega||^2}$ by \cite{Wang}, it follows that
\begin{eqnarray*}
|\omega| \Delta |\omega| - \Ric(\omega,\omega) = |\nabla \omega|^2
-|\nabla|\omega||^2 \geq {1 \over {n-1}} |\nabla|\omega||^2 .
\end{eqnarray*}
By Lemma \ref{Leung}, we have
\begin{eqnarray*}
|\omega| \Delta |\omega| - {1 \over {n-1}} |\nabla|\omega||^2 \geq
\Ric(\omega,\omega) \geq -(n-1)|\omega|^2 -{{n-1} \over n}|A|^2 |\omega|^2 .
\end{eqnarray*}
Therefore we get
\begin{eqnarray*}
|\omega| \Delta |\omega|  + {{n-1} \over n}|A|^2 |\omega|^2 + (n-1)|\omega|^2 - {1
\over {n-1}} |\nabla|\omega||^2 \geq 0.
\end{eqnarray*}
Multiplying both sides by $\varphi^2$ as in the proof of Theorem
\ref{main1} and integrating over $M$, we have from integration by
parts that
\begin{align}
0 &\leq \int_M \varphi^2 |\omega| \Delta |\omega|  + {{n-1} \over
n}\varphi^2|A|^2 |\omega|^2 + (n-1)|\omega|^2 \varphi^2 - {1
\over {n-1}} \varphi^2 |\nabla|\omega||^2 dv   \label{form}    \\
&= -2 \int_M \varphi |\omega| \langle \nabla \varphi, \nabla
|\omega| \rangle dv - {n \over {n-1}}\int_M \varphi^2
|\nabla|\omega||^2 dv \nonumber \\ & \quad + (n-1)\int_M |\omega|^2 \varphi^2 dv + {{n-1} \over n}\int_M
|A|^2 |\omega|^2\varphi^2 dv.  \nonumber
\end{align}

On the other hand, we get the following from H$\rm \ddot{o}$lder
inequality and Sobolev inequality(Lemma \ref{MS})
\begin{align*}
\int_M |A|^2 |\omega|^2\varphi^2 dv &\leq \Big(\int_M |A|^n dv
\Big)^{2 \over
n} \Big(\int_M (\varphi|\omega|)^{2n \over {n-2}}dv \Big)^{{n-2} \over n}   \\
&\leq C_s \Big(\int_M |A|^n dv \Big)^{2 \over n} \int_M |\nabla
(\varphi
|\omega|)|^2 dv \\
&= C_s \Big(\int_M |A|^n dv \Big)^{2 \over n} \Big(\int_M
|\omega|^2|\nabla \varphi|^2 + |\varphi|^2 |\nabla|\omega||^2 + 2
\varphi |\omega| \langle \nabla \varphi, \nabla |\omega| \rangle
dv \Big).
\end{align*}
Then (\ref{form}) becomes
\begin{align}
0 &\leq -2 \int_M \varphi |\omega| \langle \nabla \varphi, \nabla
|\omega| \rangle dv - {n \over {n-1}} \int_M \varphi^2 |\nabla |\omega||^2 dv + (n-1)\int_M |\omega|^2 \varphi^2 dv  \label{form1} \\
& \quad + {{n-1} \over n}C_s \Big(\int_M |A|^n dv \Big)^{2 \over n}
\Big(\int_M |\omega|^2|\nabla \varphi|^2 + \varphi^2
|\nabla|\omega||^2 + 2 \varphi |\omega| \langle \nabla \varphi, \nabla |\omega|
\rangle dv \Big)\nonumber .
\end{align}
Applying the inequality (\ref{thm : CL}), we have
\begin{eqnarray*}
\frac{(n-1)^2}{4} \leq \lambda_1 (M) \leq \frac{\int_M |\nabla (\varphi |\omega|)|^2 dv}{\int_M \varphi^2 |\omega|^2 dv} .
\end{eqnarray*}
Thus
\begin{eqnarray*}
\int_M |\omega|^2 \varphi^2  dv \leq \frac{4}{(n-1)^2} \Big(\int_M |\omega|^2|\nabla \varphi|^2 + \varphi^2
|\nabla|\omega||^2 + 2 \varphi |\omega| \langle \nabla \varphi, \nabla |\omega|
\rangle dv \Big) .
\end{eqnarray*}
Using Schwarz inequality for $\varepsilon > 0$, we have
\begin{eqnarray*}
2\Big|\int_M \varphi |\omega| \langle \nabla \varphi, \nabla
|\omega| \rangle dv\Big| \leq {\varepsilon \over 2} \int_M
\varphi^2 |\nabla|\omega||^2 dv +{2 \over \varepsilon} \int_M
|\omega|^2|\nabla \varphi|^2 dv .
\end{eqnarray*}
Therefore it follows from the inequality (\ref{form1})
\begin{align*}
&\Bigg[{n-4  \over {n-1}} - {n-1 \over n}C_s \Big(\int_M |A|^n dv
\Big)^{2 \over n} \\
&- {\varepsilon \over 2}\Big\{1 + \frac{4}{(n-1)^2} + {n-1 \over n}C_s
\Big(\int_M |A|^n dv \Big)^{2 \over n}\Big\}\Bigg] \int_M \varphi^2
|\nabla|\omega||^2 dv \\
 &\leq \Bigg[ \frac{4}{n-1} + {n-1 \over n}C_s \Big(\int_M |A|^n dv \Big)^{2 \over n} \\
 &\quad +{2 \over \varepsilon}\Big\{1 + \frac{4}{(n-1)^2} +
{n-1 \over n}\Big(\int_M |A|^n dv \Big)^{2 \over n} \Big\}   \Bigg] \int_M
|\omega|^2|\nabla \varphi|^2 dv.
\end{align*}
Since $\displaystyle{\Big(\int_M |A|^n dv \Big)^{2 \over n} < {n(n-4) \over
(n-1)^2} \frac{1}{C_s}}$ by assumption, choosing $\varepsilon
> 0$ sufficiently small and letting $R \rightarrow \infty$, we
obtain $\nabla |\omega| \equiv 0$, i.e., $|\omega|$ is constant.
However, since $\displaystyle{\int_M |\omega|^2 dv < \infty}$ and
the volume of $M$ is infinite (\cite{Anderson} and \cite{Wei}.), we get $\omega \equiv 0$. \qed

\bigskip

In the proof of Theorem \ref{main1}, if we do not use the fact that $\lambda_1 (M) \geq \frac{(n-1)^2}{4}$ and assume that
\begin{eqnarray*}
\Big(\int_M |A|^n dv \Big)^{2 \over n} < \frac{(n\lambda_1(M) - (n-1)^2)n}{(n-1)^2C_s \lambda_1 (M)}
\end{eqnarray*}
for  $\lambda_1 (M) > \frac{(n-1)^2}{n}$, one can see that there exist no nontrivial $L^2$ harmonic 1-forms on $M^n$ for $n \geq 3$ by using the same argument as in the above proof. More precisely, we have
\begin{thm}
Let $M$ be an $n$-dimensional complete immersed minimal submanifold in $\mathbb{H}^{n+p}$, $n \geq 3$. Assume that $\lambda_1 (M) > \frac{(n-1)^2}{n}$ and the total scalar curvature satisfies
\begin{eqnarray*}
\Big(\int_M |A|^n dv \Big)^{2 \over n} < \frac{(n\lambda_1(M) - (n-1)^2)n}{(n-1)^2C_s \lambda_1 (M)} .
\end{eqnarray*}
Then there are no nontrivial $L^2$ harmonic 1-forms on $M$.
\end{thm}

\end{document}